\documentstyle[12pt,amssymb,amsfonts,amsbsy,amsmath]{article}
\newtheorem{dfn}{Definition}
\newtheorem{thm}{Theorem}
\newtheorem{prp}{Proposition}

\begin{document}
\begin{center}
\bigskip
\large
ON THE GALOIS GROUP SOME FUCHSIAN SYSTEMS
\bigskip

\bigskip
\vspace*{0.3cm}  Ala Avoyan 
\bigskip
\end{center}

\vspace*{0.25cm} \noindent {\small {\bf Abstract}. {The aim of
this paper is to give a new result of the differential Galois
theory of linear ordinary differential equations. In particular,
we compute differential Galois group for special type non-resonant
Fuchsian system.}

\bigskip
\vspace*{0.25cm}\noindent{\small {\bf Keywords and phrases}:
{Galois group, Fuchsian system, monodromy representation.} }

\vspace*{0.25cm}\noindent{\small {\bf AMS subject classification
(2000):} {57R45, 12F10.}}

\vspace*{0.3cm}

\section{Introduction}\label{S1}

The Galois theory of linear ordinary differential equations was
created by Picard and Vessiot at the end of the nineteenth century
and it is a theory analogous to the classical Galois theory of
polynomials. It is known, that if given the differential field $K$
of coefficients and the extension field of $L$ of $K$ generated by
the solutions, then the essential information about the solutions
is contained in a group of hidden symmetries of the equation: the
Galois group. This group is a linear algebraic group and, as in
the classical theory, a correspondence between subfields and
subgroups is satisfied (see \cite{magid}).

 The first rigorous
proofs of nonsolvability of differential equations in a finite
form (in terms of quadratures and elementary functions) were known
before the middle of the XIXth century.  The Liouville theory does
not imply that a "simple equation" necessarily has a "simple
solution". For example, the Bessel equation
$$
y^{\prime\prime} + xy^{\prime} + (x^2-\nu^2)y = 0,
$$
which is integrable in a finite form for $\nu=\frac{2n +1 }{2},
n\in Z,$ has the solution
$$
J_{n+\frac{1}{2}}(x)=(-1)^{n+\frac{1}{2}}\frac{1}{\sqrt{\pi}}\frac{d^n}{d(x^2)^n}\frac{sinx}{x}
$$
where $n = 0, 1, 2,....$ In addition to the group analysis
mentioned above, the differential operators factorization method,
applied in combination with Kummer-Liouville and Darboux
transformations, is also an effective way of integrating 
ordinary differential equations.

Let $D = \{a_1,..., a_n\}$ be a finite set of points on the
Riemann sphere $CP^1$, and let $z$ be a parameter on $CP^1$. It is
assumed that $z=\infty$ is not among the points of $D$. Consider
Fuchsian systems of first-order linear differential equations with
the set of singularities of D, i.e., systems of linear
differential equations with first-order poles at points from $D$,
which have the form
\begin{equation}\label{fuchs}
df=\omega f
\end{equation}
where $f$ is a column vector of $p$ components, and the
coefficient matrix $\omega$ has the form
$$
\omega=\sum_{i=1}^n\frac{B_i}{z-a_i}dz
$$
with constant $p\times p $ complex matrices $B_i$ satisfying the
condition $\sum_{i=1}^n B_i=0.$  Let
\begin{equation}\label{repr}
\chi:\pi_1(CP^1-D,z_0)\rightarrow GL(n,C)
\end{equation}
be the representation of the fundamental group. The
Riemann-Hilbert problem formulated as follows: realize the
representation $\chi$ as the monodromy representation of some
Fuchsian system of the form (\ref{fuchs}) (see \cite{bolibrukh}).

According to Lappo-Danilevsky, it is possible to analytically
express coefficients of a Fuchs type systems by the monodromy
matrices, provided matrices satisfy certain conditions.
Lappo-Danilevsky showed that if the monodromy matrices
$M_1,...,M_m$ are close to $\textbf{1}$, then coefficients $A_j$
of the system of differential equations of the Fuchs type
$\frac{df}{dz}=\left(\sum_{j=1}^m\frac{A_j}{z-s_j}\right)f$ are
expressed by the singular points $s_j$ and monodromy matrices
$M_j$ via noncommutative power series $ A_j=\frac1{2\pi i}{\tilde
M}_j+\sum_{1\leq k,l\leq n}\xi_{kl}(s){\tilde M}_k{\tilde
M}_l+\cdots, $ where $\xi_{kl}$ is a function depending on the
singular points which can be given explicitly from $s\in S$, and
${\tilde M}_j=M_j-\textbf{1}$. Algebraic version of the
Riemann-Hilbert monodromy problem is known in the differential
Galois theory under the name of {\it inverse problem} [4].

The problem is formulated as follows: let $k$ be a differential
field with the field of constants $C$ and
$D(y)=y^{(n)}+a_1y^{(n-1)}+\cdots+a_{n-1}y'+a_ny$ be a
differential operator with coefficients in $k$. To the operator
$D$ one assigns the so-called {\it Picard-Vessiot field} $K$,
whose automorphism group $G$ is the Galois group of the equation
$D(y)=0$, isomorphic to some subgroup of $GL_n(C)$. More precisely
the inverse problem is:

\textit{Given a group $G$, find an extension of $k$ with Galois
group $G$.}

Constructive character of this problem will become clear if one
recalls that a system of linear differential equations is solvable
in the class of Liouville functions if and only if the identity
component of the Galois group of the equation is a solvable group
(see [5]).

Differential Galois theory is also a basis for establishing
integrability of a function in the class of elementary functions,
and after the elegant works by Kovacic, Davenport and Singer the
algorithmic approach became an alternative to the Slagle's
heuristic approach. Similar ideas appear in the theory of
integrable systems (see [5]).

Relationship between the Galois group of the regular system and
its monodromy group is expressed by the following result:

\textit{   The differential Galois group of the regular system $df=\omega f$
is the Zariski closure of its monodromy group} (see [5]).

 We begin with some basic definitions.

\begin{dfn} 1) A {\em differential ring} $(R,\Delta)$ is a
ring $R$ with a  set $\Delta = \{\partial_1, \ldots ,
\partial_m\}$ of maps ({\em derivations}) $\partial_i: R
\rightarrow R$, such that
\begin{enumerate}
\item $\partial_i(a+b) = \partial_i(a)+\partial_i(b), \ \ \partial_i(ab) = \partial_i(a)b+a\partial_i(b)$ for all $a,b\in R$, and
\item $\partial_i\partial_j = \partial_j \partial_i$ for all $i,j$.
\end{enumerate}
\noindent 2) The ring $C_R = \{ c \in R \ | \ \partial(c) = 0 \
\forall \ \partial\in \Delta\}$ is called the {\em ring of
constants of $R$}.

When $m = 1$, we say $R$ is an {\em  ordinary differential ring} $
(R,\partial)$. We frequently use the notation $a'$ to denote
$\partial(a)$ for $a \in R$. A differential ring that is also a
field is called a {\em differential field}.  If $k$ is a
differential field, then $C_k$ is also a field.
\end{dfn}

\textbf{Example.} 1)  $(C^{\infty}(R^m), \Delta =
\{\frac{\partial}{\partial x_1}, \ldots , \frac{\partial}{\partial
x_m}\}$)
 = infinitely differentiable functions on $R^m$.

2)  $(C(x_1, \ldots , x_m),  \Delta = \{\frac{\partial}{\partial
x_1}, \ldots , \frac{\partial}{\partial x_m}\}$) = field of
rational functions.

3) $(C[[x]], \frac{\partial}{\partial x})$  = ring of formal power
series $C((x)) =  \mbox{ quotient field of }C[[x]]
=C[[x]][\frac{1}{x}]$

 4) $(C\{\{x\}\}, \frac{\partial}{\partial
x})$  = ring of germs of convergent series $C(\{x\}) =  \mbox{
quotient field of }C\{\{x\}\} =C\{\{x\}\}[\frac{1}{x}]$

 5)
$(M_{O}, \Delta = \{\frac{\partial}{\partial x_1}, \ldots ,
\frac{\partial}{\partial x_m}\}$) = field of functions meromorphic
on $O^{\rm open, connected} \subset C^m$}

\section{Main Results}\label{S2}

The following result (see [5]) shows that many examples reduce to
Example 5) above:

\begin{thm} Any differential field $k$, finitely
generated over $Q$, is isomorphic to a differential subfield of
some $M_O$.
\end{thm}
We wish to consider and compare three different versions of the
notion of a linear differential equation.

\begin{prp} Let $(k,\partial)$ be a differential field.
Then the following three notions are equivalent
\begin{enumerate}
 \item A {\em scalar linear differential equation} is an equation of the form
 \[L(y) = a_ny^{(n)} + \ldots + a_0y = 0 , \ a_i \in k.\]
 \item A {\em matrix linear differential equation} is an equation of the form
 \[Y' = AY, \ A \in GL_n(k)\]
 where $GL_n(k)$ denotes the ring of $n\times n$ matrices with entries in $k$.
 \item A {\em differential module} of dimension $n$ is an $n$-dimensional $k$-vector space $M$ with a map  $\partial:M\rightarrow M$ satisfying
\[\partial(fm) = f'm+f\partial m \mbox{ for all } f \in k, m \in M.\]
\end{enumerate}
\end{prp}
Let $L(y) = y^{(n)} + a_{n-1} y^{(n-1)} + \ldots + a_0 y = 0$. If
we let
 $y_1 = y, y_2 = y', \ldots y_n = y^{(n-1)}$, then we have
\[\left(\begin{array}{c} y_1\\y_2\\ \vdots \\ y_{n-1}\\y_n \end{array}\right)' =
 \left(\begin{array}{ccccc}0&1&0&\ldots& 0\\ 0& 0 & 1& \ldots & 0 \\ \vdots &\vdots &\vdots &\vdots &\vdots
 \\ 0&0&0&\ldots& 1\\-a_0& -a_1& -a_2 & \ldots & -a_{n-1}\end{array}\right) \left(\begin{array}{c} y_1\\y_2\\ \vdots
 \\ y_{n-1}\\y_n \end{array}\right) \]
We can write the last equation as $Y' = A_LY$ and refer $A_L$
as the companion matrix of the scalar equation and the matrix
equation as the companion equation.  Clearly any solution of the
scalar equation yields a solution of the companion equation and
{\em vice versa}.

Given $Y' = AY, \ A \in GL_n(k)$,  we construct a differential
module in the following way:  Let $M = k^n, \ e_1, \ldots, e_n$
the usual basis.  Define $\partial e_i = -\sum_j a_{j,i} e_j, $
{\em i.e.,} $\partial e = -A^t e$.  Note that if $m =
\sum_if_ie_i$ then $\partial m = \sum_i(f_i' -
\sum_ja_{i,j}f_j)e_i$. In particular, we have that $\partial m =
0$ if and only if
\[\left(\begin{array}{c} f_1\\ \vdots \\ f_n\end{array}\right)' = A  \left(\begin{array}{c} f_1\\ \vdots \\ f_n\end{array}\right)\]
It is this latter fact that motivates the seemingly strange
definition of this differential module, which we denote by $(M_A,
\partial)$.

Conversely, given a differential module $(M, \partial)$,
\textit{select} a basis $e = (e_1, \ldots ,e_n)$.  Define $A_M \in
GL_n(k)$ by $\partial e_i = \sum_j a_{j,i}e_j$.  This yields a
matrix equation $Y' = AY$.  If $\bar{e} = (\bar{e}_1,  \ldots
,\bar{e}_n)$ is another basis, we get another equation $Y' =
\bar{A}Y$.  If $f$ and $\bar{f}$ are vectors with respect to these
two bases and $f = B\bar{f}, \ B \in GL_n(k)$, then
\[\bar{A} = B^{-1}AB - B^{-1}B' \ . \]

\begin{thm} Let the system (1) be non-resonant. Then the
differential Galois group is generated by $e^{2\pi i B_j},$
$j=1,...,m.$
\end{thm}

Firstly, we note that in non-resonant case the differential Galois
group of equation coincides with closure of the subgroup of
$GL_n(C)$ generated by monodromy matrices of the system. Next, use
expression
\begin{equation}\label{i5}
\Phi_j(\widetilde{z})=U_j(z)(z-s_j)^{\Lambda_j}(\widetilde{z}-s_j)^{E_j},
\end{equation}
of the solution to system (1) and the fact, that a) the sums
$\rho_j^i+\varphi_j^i$ are the eigenvalues of the matrix-residue
and b) equal eigenvalues of $E_j$ occupy consecutive position on
the diagonal and that the matrix $E_j$ is block-diagonal, with
diagonal blocks of sizes equal to their multiplicities. Here
$\widetilde{z}$ denotes the coordinates  on the universal
covering, $U_j$ is holomorphic invertible at $s_j,$
$\Lambda_j=diag(\alpha_j^1,...,\alpha_j^n),\alpha_j^1\geq...\geq\alpha_j^n,$
is the diagonal matrix, with integer entries and $E_j$ are upper
triangular matrices.
 Hence, if the eigenvalues of $B_j$ are non-resonant, then to
 equal eigenvalues of $E_j$ correspond equal eigenvalues of
 $A_j,$ the matrices $A_j$ and $E_j$ commute and
 $B_j=U_j(0)(A_j+E_j)U_j(0)^{-1}$ (see [1],[2],[3]). One has $M_j=G_j^{-1}e^{2\pi i
 E_j}G_j.$ Consequently, for the Jordan normal form $JNF(M_j)$ of
 $M_j$ one has $JNF(M_j)=JNF(E_j)=JNF(A_j+E_j)=JNF(B_j).$ The
 theorem is proved.

\textbf{Acknowledgements.} The author is grateful to
the participants of Tbilisi State University seminar "Elliptic systems on Riemann surfaces"
for their useful comments and suggestions.

\end{document}